\documentclass[12pt,fleqn]{article}
\usepackage{amsfonts}
\usepackage{epsfig}
\usepackage{float}
\usepackage{dsfont}
\usepackage{amssymb}
\usepackage{mathrsfs}
\usepackage{enumerate}
\usepackage{hhline}
\usepackage{amsmath}
\usepackage{lineno}
\usepackage{indentfirst}

\usepackage{graphicx}
\newtheorem{proposition}{Proposition}[section]

\newtheorem{lemma}[proposition]{Lemma}

\newtheorem{defi}[proposition]{Definition}
\newtheorem{fact}[proposition]{Fact}
\newtheorem{claim}[proposition]{Claim}

\newtheorem{theorem}[proposition]{Theorem}

\newtheorem{problem}[proposition]{Problem}

\newtheorem{corollary}[proposition]{Corollary}

\def\p{\noindent{\bf Proof. }}
\def\q{\hspace*{\fill}$\Box$\medskip}
\usepackage{color}
\usepackage{cite}

\begin{document}
\title{ {\bf Bipartite Ramsey numbers of large cycles}}
\author{Shaoqiang Liu \thanks{ College of Mathematics and Econometrics, Hunan University, Changsha 410082, P.R. China. Email: hylsq15@sina. com.}
\and Yuejian Peng \thanks{Corresponding author. Institute of Mathematics, Hunan University, Changsha 410082, P.R. China. Email: ypeng1@hnu.edu.cn.
Partially  supported by National Natural Science Foundation of China (No. 11671124).} }

\date{}
\maketitle
\vskip 0.1in
\begin{abstract}

For an integer $r\geq 2$ and bipartite graphs $H_i$, where $1\leq i\leq r$, the bipartite Ramsey number $br(H_1,H_2,\ldots,H_r)$ is the minimum integer $N$ such that any $r$-edge
coloring of the complete bipartite graph $K_{N,N}$ contains a monochromatic subgraph isomorphic to $H_i$ in color $i$ for some $i$, $1\leq i\leq r$. We show that for
$\alpha_1,\alpha_2>0$,
$br(C_{2\lfloor \alpha_1 n\rfloor},C_{2\lfloor \alpha_2 n\rfloor})=(\alpha_1+\alpha_2+o(1))n$. We also show that if $r\geq 3, \alpha_1,\alpha_2>0,  \alpha_{j+2}\geq
[(j+2)!-1]\sum^{j+1}_{i=1}
\alpha_i$ for $j=1,2,\ldots,r-2$, then $br(C_{2\lfloor \alpha_1 n\rfloor},C_{2\lfloor \alpha_2 n\rfloor},\ldots,C_{2\lfloor \alpha_r n\rfloor})=(\sum^r_{j=1} \alpha_j+o(1))n.$
For
$\xi>0$ and sufficiently large $n$, let $G$ be a bipartite graph with bipartition $\{V_1,V_2\}$, $|V_1|=|V_2|=N$, where $N=(2+8\xi)n$. We prove that if
$\delta(G)>(\frac{7}{8}+9\xi)N$, then any $2$-edge coloring of $G$ contains a monochromatic copy of $C_{2n}$.
\vskip2mm
\noindent\textbf{Key words:} Cycle; Ramsey number;  Bipartite ramsey number; Regularity lemma

\end{abstract}

\section{Introduction}

\baselineskip 14pt

Let $r\geq 2$ be an integer and  $H_1,\ldots,H_r, H$ be given graphs. The {\em Ramsey number} $R(H_1,\ldots,H_r)$ is the minimum integer $N$ such that any
$r$-edge coloring of the complete graph $K_N$ contains a monochromatic subgraph in color $i$ isomorphic to $H_i$ for some $i$, $1\leq i\leq r$.
$R(H_1,\ldots,H_r)$ is simplified by $R_r(H)$ if $H_i=H$ for every $1\leq i \leq r$. Let $B_1,\ldots,B_r, B$ be bipartite graphs. The {\em bipartite Ramsey number} $br(B_1,\ldots, B_r)$ is
the minimum integer $N$ such that any
$r$-edge coloring of the complete bipartite graph $K_{N, N}$ contains a monochromatic subgraph in color $i$ isomorphic to $B_i$ for some $i$, $1\leq i\leq r$. If $B_i=B$ for every
$1\leq i\leq r$, then it is simplified by $br_r(B)$. We refer the reader to the book by Graham, Rothschild and Spencer \cite{GRS} for an overview and a survey by Conlon, Fox and Sudakov \cite{CFS} for recent developments in Ramsey theory. In this paper, we focus on bipartite Ramsey number of cycles.

Let $C_n$ be a cycle of length $n$. The Ramsey number $R(C_m,C_n)$ has been completely determined by several authors, including Bondy and Erd\H{o}s~\cite{JB73}, Faudree and Schelp~\cite{RF74}, and
Rosta~\cite{VR731,VR732}. Bondy and Erd\H{o}s~\cite{JB73} conjectured that if $k\geq 2$ and $n\geq 3$ is odd, then $R_{k}(C_n)=2^{k-1}(n-1)+1$. {\L}uczak~\cite{TL99} proved that for
$n$ odd, $R_{3}(C_n)=4n+o(n)$.
Recently, Jenssen and Skokan~\cite{MJ16} resolved the above conjecture of  Bondy and Erd\H{o}s~\cite{JB73} for large $n$.  Benevides and Skokan~\cite{FS09} proved that there exists
$n_1$
such that for every even $n\geq n_1$, $R_3(C_{n})=2n$. S\'ark\"ozy~\cite{GS16} showed that, for every $r\geq2$ and even $n$, $R_r(C_{n})\leq (r-\frac{r}{16r^3+1})n+o(n)$.  Davies et
al.~\cite{ED17} showed that $R_r(C_{n})\leq (r-\frac{1}{4})n+o(n)$ for $r\geq4$. Recently, Knierim and Su ~\cite{CK18} improved the upper bound on the multicolour Ramsey nunbers of
even cycles to $R_r(C_{n})\leq (r-\frac{1}{2})n+o(n)$ for $r\geq4$.

The bipartite Ramsey numbers for cycles has also been studied actively.  A lower bound
$$b_{r}(C_{2n_1},C_{2n_2},\\ \ldots,C_{2n_r})\geq \sum_{i=1}^{r}n_i-r+1$$ is implied by the following
construction: Let $N= \sum_{i=1}^{r}n_i-r$. Take the complete bipartite graph $K_{N, N}$ with the vertex set $U=\{u_1,u_2,\ldots,u_N\},
V=\{v_1,v_2,\ldots,v_N\}$. Partition it into $r$ color classes $H_k(k=1,2,3,\ldots,r)$, in which $u_sv_t$ is an edge in $H_k$ if and only if $u_s\in U$ and
$\sum_{i=1}^{k-1}n_i-k+2\leq t\leq \sum_{i=1}^{k}n_i-k$. Clearly, $H_k$ dose not contain $C_{2n_k}$
as a subgraph. This gives the lower bound.

Zhang et al.~\cite{RZ11,RZ13} determined $br(C_{2n},C_4)=n+1$, $br(C_{2n},C_6)=n+2$ for $n\geq 4$. Goddard
et
al.~\cite{WG08} determined $br(C_4,C_4,C_4)=11$. Joubert~\cite{EJ17} showed that $$br(C_{2t_1},C_{2t_2},\ldots,C_{2t_k})\leq k(t_1+t_2+\cdots+t_k-k+1),$$ where $t_i$ is an integer and
$2\leq t_i\leq 4$, for all $1\leq i\leq k$. Recently, Shen, Lin and Liu~\cite{LS18} showed that $br(C_{2n},C_{2n})=2n+o(n)$. Let $t_i\geq 3$ be positive integers for all $1\leq i\leq
k$,
Hattingh and Joubert~\cite{JH18} proved that $br(C_{2t_1},C_{2t_1},\ldots,C_{2t_k})\leq k(2t-4k+1)$, where $t=t_1+t_2+\cdots+t_k$. Buci\'c, Letzter and Sudokov~\cite{MB18} showed that
$br_3(C_{2n})=(3+o(1))n$.  The authors~\cite{LP18} showed that $br_r(C_{2n})\leq r(1+\sqrt{1-\frac{2}{r}}) n+o(n)$ for $r\geq 2$. In this paper, we show the following result.

\begin{theorem}\label{1-1}
For $\xi>0$ and $\alpha_1,\alpha_2>0$, there exists $n_0$ such that for any $n\geq n_0$ the following holds: If $N\geq (\alpha_1+\alpha_2+8\xi)n$, then any $2$-edge
coloring of $K_{N,N}$ contains a monochromatic copy of  $C_{2\lfloor \alpha_1 n\rfloor}$ in color $1$ or $C_{2\lfloor \alpha_2 n\rfloor}$ in color $2$.
\end{theorem}

The following result is implied by Theorem~\ref{1-1}.

\begin{corollary}\label{1-2} For $\alpha_1,\alpha_2>0$ , we have
$$br(C_{2\lfloor \alpha_1 n\rfloor},C_{2\lfloor \alpha_2 n\rfloor})=(\alpha_1+\alpha_2+o(1))n.$$
\end{corollary}

We also obtain the following result on bipartite Ramsey numbers of large cycles in multicolorings.

\begin{theorem}\label{1-3}
For $r\geq 3$, $\alpha_1,\alpha_2>0$, $\alpha_{j+2}\geq [(j+2)!-1]\sum^{j+1}_{i=1} \alpha_i, j=1,2,\ldots,r-2$ and $\xi>0$, there exists $n_0$ such that for any $n\geq n_0$ the
following holds: If $N\geq (\sum^r_{i=1} \alpha_i+(r+2)r!\xi)n$, then any $r$-edge coloring of $K_{N,N}$ contains a
monochromatic copy of $C_{2\lfloor \alpha_i n\rfloor}$ for some $1\leq i\leq r$.
\end{theorem}

The following result is implied by Theorem~\ref{1-3}.

\begin{corollary}\label{1-4}  For $r\geq 3, \alpha_1,\alpha_2>0, \alpha_{j+2}\geq [(j+2)!-1]\sum^{j+1}_{i=1} \alpha_i, j=1,2,\ldots,r-2$, we have
$$br(C_{2\lfloor \alpha_1 n\rfloor},C_{2\lfloor \alpha_2 n\rfloor},\ldots,C_{2\lfloor \alpha_r n\rfloor})=(\sum^r_{j=1} \alpha_j+o(1))n.$$
\end{corollary}

For a graph $G$, denote by $V(G)$ the set of vertices of $G$ and $\delta(G)$ the minimum degree. Motivated by a new class of Ramsey-Tur\'an type problems raised by
Schelp~\cite{RS12}. The intent is to find the smallest positive constant $0<c<1$, such that if $\delta(G)\geq c|V(G)|$, then any $r$-edge-coloring of $G$ contains a monochromatic subgraph $H$. If
$\alpha_1=\alpha_2=1$, then Theorem~\ref{1-1} can be strengthened to the following result.

\begin{theorem}\label{1-5}
For  $\xi>0$, there exists $n_0$ such that for any $n\geq n_0$ the following holds: Let $G$ be a bipartite graph with bipartition $\{V_1,V_2\}$, $|V_1|=|V_2|=N=(2+8\xi)n$,
and $\delta(G)>(\frac{7}{8}+9\xi)N$. Then any $2$-edge coloring of $G$ contains a monochromatic copy of $C_{2n}$.
\end{theorem}

We will give the proofs of Theorems~\ref{1-1},~\ref{1-3} and~\ref{1-5} in Sections 2, 3 and 4 respectively. The basic idea of the proof of  Theorems~\ref{1-1} and ~\ref{1-5}
follows  from the work of Figaj and {\L}uczak~\cite{FT07}: Apply the bipartite form of the regularity Lemma to $K_{N,N}$, consider the auxiliary graph by viewing each part in the regularity
partition as a vertex, show the existence of a ``fat" connected matching and expand it to a ``long" cycle.  The crucial part is to show the existence of a ``fat" connected matching
in the auxiliary graph. In~\cite{FT07}, the host graph is a complete graph and the auxiliary graph is an ``almost" complete graph. In our case, the host graph is changed to a
complete bipartite graph (or a
bipartite graph with ``large"  minimum degree) and the auxiliary graph is an ``almost" complete  bipartite graph (or a
bipartite graph with ``large"  minimum degree).

\section{2-color bipartite Ramsey number of cycles}

For a graph $G$, let $E(G)$ be the set of its edges and $e(G)=|E(G)|$. Denote the degree of a vertex $v$ in a graph $G$ by $d_G(v)$ and the neighbors of $v$ in a graph $G$ by $N(v)$.
We
shall denote by $G[H]$ the induced subgraph of $G$ by $H\subseteq V(G)$. A set $M$ of pairwise disjoint edges of a graph $G$ is called a matching. $M$ saturates
a vertex $u$ or $u$ is $M$-\emph{saturated} if there is an edge in $M$ incident with $u$; otherwise, $u$ is $M$-\emph{unsaturated}. A matching $M$ is connected in $G$ if all edges of
$M$ are in the same component of $G$. If the edges of $G$ are $r$-colored, then let $G_{i}$ denote the spanning subgraph of $G$ with all edges colored by $i$.

\subsection{The Regularity Lemma}

Let $A,B$ be disjoint subsets of $V(G)$. Let $e(A,B)=e_G(A,B)$ denote the number of edges $\{u,v\}$ with $u\in A$ and $v\in B$. The ratio
$$
d(A,B)=\frac{e(A,B)}{|A||B|}
$$
is called the edge density of $(A,B)$. Clearly, $0\leq d(A,B)\leq 1$.

For $\epsilon>0$, a disjoint pair $(A,B)$ is called $\epsilon$-regular if
$
 |d(A,B)-d(A^{\prime},B^{\prime})|\leq\epsilon
$ for any $A^{\prime}\subseteq A$ and $B^{\prime}\subseteq B$ with $|A^{\prime}|>\epsilon |A|$ and $|B^{\prime}|>\epsilon|B|$.
\begin{fact}\label{2-0}
Let $(A,B)$ be an $\epsilon$-regular pair with density $d$ and $B^{\prime}\subseteq B$ with $|B^{\prime}|>\epsilon|B|$, then all but at most $\epsilon|A|$ vertices $v\in |A|$ satisfy
$|N(v)\bigcap B^{\prime}|>(d-\epsilon)|B^{\prime}|$.
\end{fact}

Benevides and Skokan~\cite{FS09} gave the next lemma which is a slightly stronger version of Claim 3 in~\cite{TL99}.

\begin{lemma}[\cite{FS09}]\label{2-1}
For every $0<\beta\leq 1$ there exists $n_0$ such that for every $n>n_0$ the following holds: Let $G$ be a bipartite
graph with bipartition $V(G)=V_1\bigcup V_2$ such that $|V_1|,|V_2|=n$. Furthermore let the pair $(V_1,V_2)$ be
$\epsilon$-regular with density at least $\frac{\beta}{4}$ for some $\epsilon$ satisfying $0<\epsilon<\frac{\beta}{100}$.
Then for every $l$, $1\leq l\leq n-\frac{5\epsilon n}{\beta}$, and for every pair of vertices $v^{\prime}\in V_1, v^{\prime\prime}\in V_2$
satisfying $d_G(v^{\prime}), d_G(v^{\prime\prime})\geq \frac{\beta n}{5}$, $G$ contains a path of length $2l+1$ connecting $v^{\prime}$
and $v^{\prime\prime}$.
\end{lemma}

 We will use the following multicolored regularity lemma for bipartite graphs.

\begin{lemma}[\cite{HL09,QL15}]\emph{(}Multicolored regularity lemma for bipartite graphs\emph{)}\label{2-2}
For any real $\epsilon>0$ and integers $K_0\geq1$ and $r\geq 1$, there exists $M=M(\epsilon,K_0)$ such that if the edges of a bipartite graph
$G=G(W^{(1)},W^{(2)})$ with $|W^{(1)}|=|W^{(2)}|\geq M$ are $r$-colored, then there exists a partition $\{W^{(s)}_1,\ldots,
W^{(s)}_k\}$ for each $W^{(s)}(s=1,2)$, $K_0\leq k \leq M$, such that
\begin{enumerate}[(1)]
  \item $||W^{(s)}_i|-|W^{(s)}_j||\leq1$ for each $s$;
  \item All but at most $\epsilon rk^2$ pairs $(V^{(1)}_i,V^{(2)}_j),1\leq i,j\leq k$, are $\epsilon$-regular for each monochromatic graph $G_1,G_2,\ldots,G_r$.
\end{enumerate}
\end{lemma}

\subsection{Proof of Theorem~\ref{1-1}}

 We may assume that $0<\xi\leq 10^{-3}$. Let $0<\epsilon \leq \frac{1}{2}\xi^7$. Let $n$ be a sufficiently large positive integer. Consider a $2$-edge coloring of $K_{N,N}$, where
 $N=(\alpha_1+\alpha_2+\xi)n$. Lemma~\ref{2-2} guarantees that there exists a partition $\{W^{(s)}_1,\ldots,W^{(s)}_k\}$ of $W^{(s)}(s=1,2)$, $N_0\leq k \leq M$, where $N_0$ is
 from Lemma~\ref{2-3}, such that

\begin{enumerate}[(1)]
  \item $\big||W^{(s)}_i|-|W^{(s)}_j|\big|\leq1$ for each $s$;
  \item All but at most $2\epsilon k^2$ pairs $(V^{(1)}_i,V^{(2)}_j),1\leq i,j\leq k$, are $\epsilon$-regular for each monochromatic graph $G_1,G_2$.
\end{enumerate}

Now we define the reduced bipartite graph $H$ in the following way: the vertex set of $H$ is $\{v^{(s)}_{1},v^{(s)}_{2},\ldots,v^{(s)}_{k}\}$, where $s=1$ or $2$,
and the edge set is defined as
\begin{center}
$E(H)=\{v^{(1)}_{i}v^{(2)}_{j}:(V^{(1)}_i,V^{(2)}_j)$ is $\epsilon$-regular for each $G_1,G_2\}.$
\end{center}
\noindent Note that by $(2)$,
$$e(H)\geq (1-2\epsilon ) k^2.$$

We define a $2$-edge coloring of $H=H_1\bigcup H_2$ in the following way: for $f\in \{1,2\}$, we put the edge $v^{(1)}_{i}v^{(2)}_{j}$ into  $H_{f}$ (lexicographically first) if
$e_{G_{f}}(V^{(1)}_i,V^{(2)}_j)\geq (\frac{1}{2}-\epsilon)|V^{(1)}_i||V^{(2)}_j|$.

The following crucial Lemma guarantees the existence of a large connected monochromatic matching in the $2$-edge coloring of $H=H_1\bigcup H_2$.

\begin{lemma}\label{2-3}
For $0<\alpha_2\leq \alpha_1\leq 1$ and $0<\xi<10^{-3}$, there exist
$N_0$, such that for each bipartite graph with bipartition $\{V_1,V_2\}$, $|V_1|=|V_2|=(\alpha_1+\alpha_2+8\xi)k^{\prime}\geq N_0$ and
$|E(G)|\geq (1-\xi^7)|V_1|^2$, the following holds. For every $2$-edge coloring of $G=G_1\bigcup G_2$ there exists a color $i$, $i\in\{1,2\}$,
such that some component of $G_i$ contains a matching saturating at least $(2\alpha_i+0.1\xi)k^{\prime}$
vertices.
\end{lemma}

The proof of Lemma~\ref{2-3} will be given in subsection 2.3.

Apply Lemma~\ref{2-3} with $k=(\alpha_1+\alpha_2+8\xi)k^{\prime}$ to the $2$-edge coloring of $H=H_1\bigcup H_2$. Thus there exists some $f\in\{1,2\}$, such that one of the
components
of $H_f$ contains a matching $M^*$ saturating at least $2a=(2a_f+0.1\xi)k^{\prime}$ vertices. Let $F^*$ be a minimal
connected tree of $H_f$ containing $M^*$. Consider a closed walk $W=v^{(s)}_{i_1}v^{(3-s)}_{i_2}\cdots v^{(3-s)}_{i_t}v^{(s)}_{i_1}$ which contains all edges of $M^*$. Since $F^*$ is
a tree, so $W$ must be of even length. Applying Fact~\ref{2-0} repeatedly, we can show that there exist $v^{(q)}_{i_r}\in V^{(q)}_{i_r}$, $1\leq r\leq t+1, q=s$ or $3-s$,
$i_{t+1}=i_1$, satisfying:

\begin{enumerate}[(1)]
  \item $v^{(q)}_{i_r}$ has at least $\frac{1}{5}|V^{(q)}_{i_r}|$ neighbors in both $V^{(3-q)}_{i_{r-1}}$ and $V^{(3-q)}_{i_{r+1}}$ for each $r$;
  \item If $v^{(q)}_{i_r}v^{(3-q)}_{i_{r+1}}$ is not an edge in $M^{*}$(view an edge in the closed walk as in $M^{*}$ only when it is in $M^{*}$ and the closed walk first passes
      it), then $v^{(q)}_{i_r}v^{(3-q)}_{i_{r+1}}$ is an edge in $H_f$.
\end{enumerate}

 Let $m=\lfloor\frac{N}{k}\rfloor$.
Applying Lemma~\ref{2-1} with $\beta=1$, for every $l$, $1\leq l\leq m-5\epsilon m$, each edge $v^{(q)}_{i_r}v^{(3-q)}_{i_{r+1}}$ in $M^*$ can be enlarged to a path of length $2l+1$
connecting $u^{(q)}_{i_r}\in V^{(q)}_{i_r}$
and $u^{(3-q)}_{i_{r+1}}\in V^{(3-q)}_{i_{r+1}}$. Therefore, there exists a cycle in
$G_f$ of each even length $l^{\prime}= \sum^{a}_{j=1}2l_j+t$, where $1\leq l_j\leq m-5\epsilon m$ for $1\leq j\leq a$,  and
\begin{align*}
\sum^{a}_{j=1}2(1-5\epsilon)m+t&\geq 2a(1-5\epsilon)\frac{N(1-\epsilon)}{k}+t\\
&\geq(2a_f+0.1\xi)k^{\prime}(1-5\epsilon)(1-\epsilon)\frac{N}{k}+t\\
&\geq(2a_f+0.1\xi)(1-16\epsilon)n+t\\
&>(2a_f+0.05\xi)n+t,
\end{align*}
\noindent where the last inequality holds when $0<\epsilon \leq \frac{1}{2}\xi^7$.

Therefore, $G_f$ contains a cycle of length $2\lfloor a_f n\rfloor$ for some $f\in \{1,2\}$ and so does $G$.

\subsection{Proof of Lemma~\ref{2-3}}
We will apply the following results given by Figaj and {\L}uczak in~\cite{FT07} . The first one is a
direct consequence of Tutte's $1$-factor Theorem.

\begin{lemma}[\cite{FT07}]\label{2-4}
If a graph $G=(V,E)$ contains no matching saturating at least $\alpha$ vertices, then there exists a partition $\{S,T,U\}$
of $V$ such that
\begin{enumerate}[(1)]
\item the subgraph induced in $G$ by $T$ has maximum degree less than $\sqrt{|V|}-1$;
\item $G$ contains no edges joining $T$ and $U$;
\item $|U|+2|S|<\alpha +\sqrt{|V|}$.
\end{enumerate}
\end{lemma}

\begin{lemma}[\cite{FT07}]\label{2-5}
Let $G=(V,E)$ be a bipartite graph with bipartition $\{V_1,V_2\}$,$|V_1|\geq |V_2|$, and at least $(1-\xi)|V_1||V_2|$ edges for some $0<\xi<0.01$.
Then, there is a component in $G$ of at least $(1-3\xi)(|V_1|+ |V_2|)$ vertices containing a matching of size at least $(1-3\xi)|V_2|$.
\end{lemma}

We will show the following result.

\begin{lemma}\label{2-6}
Let $0<\alpha_2\leq \alpha_1\leq 1$ and  $0<\xi\leq10^{-3}$. Let $G=(V,E)$ be a bipartite graph with bipartition $\{V_1,V_2\}$, $|V_1|=|V_2|=(\alpha_1+\alpha_2+8\xi)k^{\prime}$ and
$|E(G)|\geq (1-\xi^2)|V_1|^2$. Then any $2$-edge coloring of $G$ leads to a monochromatic component with at least $(\alpha_1+\alpha_2+3\xi)k^{\prime}$ vertices.
\end{lemma}

\p Let us consider a $2$-edge coloring of $G$, say, red and blue. Let $F_1$ be a monochromatic component in $G$ with largest number of vertices, say, red. If $|V(F_1)|\geq
(\alpha_1+\alpha_2+3\xi)k^{\prime}$, then we are done. Hence we may assume that $|V(F_1)|<(\alpha_1+\alpha_2+3\xi)k^{\prime}$.

Clearly,  $|V(F_1)\bigcap V_s|<(\alpha_1+\alpha_2+3\xi)k^{\prime}$ for each $s=1,2$. Note
that $G$ contains a vertex of degree larger than $(\alpha_1+\alpha_2+4\xi)k^{\prime}$, so $|V(F_1)|>\frac{1}{2}(\alpha_1+\alpha_2+4\xi)k^{\prime}$.

We claim that $|V(F_1)\bigcap V_s|\geq \frac{1}{4}(\alpha_1+\alpha_2+4\xi)k^{\prime}$ for each $s=1,2$. Without loss of generality, assume that $|V(F_1)\bigcap V_1|<
\frac{1}{4}(\alpha_1+\alpha_2+4\xi)k^{\prime}$. Thus $|V_1\backslash V(F_1)|\geq \frac{3}{4}(\alpha_1+\alpha_2+9\xi)k^{\prime}$. Since
$|V(F_1)|>\frac{1}{2}(\alpha_1+\alpha_2+4\xi)k^{\prime}$, $|V(F_1)\bigcap V_2|\geq \frac{1}{4}(\alpha_1+\alpha_2+4\xi)k^{\prime}$. Clearly, $\xi^2|V_1|^2<10\xi^2|V_1\backslash
V(F_1)||V(F_1)\bigcap V_2|$. By lemma~\ref{2-5}, $G[V_1\backslash V(F_1),V(F_1)\bigcap V_2]$ contains a blue component with at least $(1-30\xi^2)(|V_1\backslash
V(F_1)|+|V(F_1)\bigcap V_2|) > (\alpha_1+\alpha_2+3\xi)k^{\prime}$ vertices, a contradiction to the choice of $F_1$.

Since $|V(F_1)|<|V_s|$, $|V_{3-s}\backslash V(F_1)|>|V(F_1)\bigcap V_s|\geq \frac{1}{4}(\alpha_1+\alpha_2+4\xi)k^{\prime}$ for each $s=1,2$.  Thus
$\xi^2|V_1|^2<64\xi^2|V(F_1)\bigcap V_{3-s}||V_s\backslash V(F_1)|.$ By lemma~\ref{2-5},
we can take $H_{1,3-s}\subseteq V(F_1)\bigcap V_{3-s}$, $H_{2,s}\subseteq V_s\backslash V(F_1)$, $|H_{1,3-s}|=|V(F_1)\bigcap V_{3-s}|-\xi k^{\prime}$, $|H_{2,s}|= |V_s\backslash
V(F_1)|-\xi k^{\prime}$ such that  $G[H_{1,3-s},H_{2,s}]$ is blue and connected. If $G[H_{2,1},H_{2,2}]$ has at least a blue edge, then $G[H_{1,2}\bigcup H_{2,1}$, $H_{1,1}\bigcup
H_{2,2}]$
is contained in a blue component with at least $2(\alpha_1+\alpha_2+6\xi)k^{\prime}$ vertices, a contraction to the choice of $F_1$. Otherwise, note that
$|H_{2,1}|+|H_{2,2}|=|V_1\backslash
V(F_1)|+|V_2\backslash V(F_1)|-2\xi k^{\prime} > (\alpha_1+\alpha_2+11\xi)k^{\prime}$ and $\xi^2|V_1|^2< \xi |H_{2,1}||H_{2,2}|$, by lemma~\ref{2-5}, $G[H_{2,1}, H_{2,2}]$ contains a
red component with at least $(1-3\xi)(|H_{2,1}|+|H_{2,2}|)\geq (1-3\xi)(\alpha_1+\alpha_2+11\xi)k^{\prime}> (\alpha_1+\alpha_2+3\xi)k^{\prime}$ vertices, a contradiction to the
choice of $F_1$. \q

Now we will give the proof of Lemma~\ref{2-3}.

\vskip3mm

\noindent\textbf{Proof of Lemma~\ref{2-3}.} For $0<\alpha_2\leq \alpha_1\leq 1$, let $G$ be a bipartite graph with bipartition $\{V_1,V_2\}$,
$|V_1|=|V_2|=(\alpha_1+\alpha_2+8\xi)k^{\prime}\geq N_0$ and $|E(G)|\geq (1-\xi^7)|V_1|^2$ for some $0<\xi\leq 10^{-3}$. Let us consider a $2$-edge coloring of $G$, say, red and
blue. Let $F_1$ denote a monochromatic component of $G$ with the largest number
of vertices and $C=V(F_1)$. Without loss of generality, assume that the edges of $F_1$ are colored by red and $F_1$ does not contain a connected matching saturating at least
$(2\alpha_1+0.1\xi)k^{\prime}$ vertices. By Lemma~\ref{2-6}, we have
$$|C|\geq(\alpha_1+\alpha_2+3\xi)k^{\prime}.$$

\textbf{Case 1} $|C|<(2\alpha_1+\alpha_2+10\xi)k^{\prime}$.

If $|C\bigcap V_1|$ or $|C\bigcap V_2|< (\alpha_1+\xi)k^{\prime}$, then, without loss of generality, assume that $|C\bigcap V_1|<(\alpha_1+\xi)k^{\prime}$. Clearly, $|V_2\bigcap
C|=|C|-|V_1\bigcap C|> (\alpha_2+2\xi)k^{\prime},|V_1\backslash C|=|V_1|-|V_1\bigcap C|> (\alpha_2+7\xi)k^{\prime}$. Since $\xi^7(\alpha_1+\alpha_2+10\xi)^2(k^{\prime})^2\leq
\xi^5|V_2\bigcap C||V_1\backslash C|$,  by lemma~\ref{2-5}, $G[V_1\backslash C,V_2\bigcap C]$ contains  a blue connected matching saturating at least $2(1-3\xi^5)\min\{|V_1\backslash
C|,|C\bigcap V_2|\}> (2\alpha_2+3\xi) k^{\prime}$ vertices.
Hence we may assume that $|V_1\bigcap C|$, $|V_2\bigcap C|\geq (\alpha_1+\xi)k^{\prime}$.

If $|V_{3-s}\backslash C|\geq (\alpha_2+\xi)k^{\prime}$ for some $s$, $s=1,2$, then, similarly, $G[V_s\bigcap C,V_{3-s}\backslash C]$ contain a blue connected matching saturating at
least $(2\alpha_2+\xi)k^{\prime}$ vertices. Hence we may assume that $|V_1\backslash C|,|V_2\backslash C|< (\alpha_2+\xi)k^{\prime}$.

Note that $|V_1\backslash C|,|V_2\backslash C|\geq5\xi k^{\prime}$. Suppose not, either $|V_1\backslash C|<5\xi k^{\prime}$ or $|V_2\backslash C|<5\xi k^{\prime}$. Without loss of
generality, assume that $|V_1\backslash C|<5\xi k^{\prime}$. Thus $|C\bigcap V_1|\geq (\alpha_1+\alpha_2+3\xi)k^{\prime}$ and hence $|C\bigcap V_2|=|C|-|C\bigcap V_1|<(\alpha_1+7\xi)
k^{\prime}$ and $|V_2\backslash C|> (\alpha_2+\xi)k^{\prime}$, a contradiction.

If $G[C\bigcap V_1, C\bigcap V_2]$ contains at least $\xi^4(k^{\prime})^2$ blue edges, then, since $\xi^7(\alpha_1+\alpha_2+10\xi)^2(k^{\prime})^2< 2\xi^5\min\{|V_1 \backslash
C||C\bigcap V_2|,|V_2 \backslash C||C\bigcap V_1|\}$, by lemma~\ref{2-5}, we can take $H^{\prime}_1\subseteq V_1\backslash C,H^{\prime}_2\subseteq C\bigcap V_2
,H^{\prime\prime}_1\subseteq V_2\backslash C,H^{\prime\prime}_2\subseteq C\bigcap V_1, |H^{\prime}_1|= |H^{\prime}_2|>|V_1\backslash C|-6\xi^5 k^{\prime},|H^{\prime\prime}_1|=
|H^{\prime\prime}_2|>|V_2\backslash C|-6\xi^5 k^{\prime}$ such that $G[H^{\prime}_1,H^{\prime}_2]$, $G[H^{\prime\prime}_1, H^{\prime\prime}_2]$ are connected and contains a perfect
matching, and  $G[H^{\prime}_2,H^{\prime\prime}_2]$ contains at least an blue edge. Thus $G[H^{\prime}_1\bigcup H^{\prime\prime}_2$, $H^{\prime}_2\bigcup H^{\prime\prime}_1]$
contains a blue connected matching saturating at least $2(H^{\prime}_1+H^{\prime\prime}_1)\geq 2(2|V_1|-12\xi^5 k^{\prime}-|C|)>2(\alpha_2+5\xi)k^{\prime}$ vertices. Hence we can
assume that $G[C\bigcap V_1, C\bigcap V_2]$ contain at most $\xi^4(k^{\prime})^2$ blue edges.

All but at most $\xi^4(k^{\prime})^2+\xi^7(\alpha_1+\alpha_2+10\xi)^2(k^{\prime})^2\leq 2\xi^4(k^{\prime})^2$ edges in $G[C\bigcap V_1, C\bigcap V_2]$ are red.  By lemma~\ref{2-5},
$G[C\bigcap V_1, C\bigcap V_2]$ contains a red connected matching saturating at least $2(1-6\xi^2)\min\{|C\bigcap V_1|$, $|C\bigcap V_2|\}> (2\alpha_1+\xi)k^{\prime}$ vertices, a
contradiction.

\textbf{Case 2} $|C|\geq(2\alpha_1+\alpha_2+10\xi)k^{\prime}$.

Since $F_1$ contains no matching saturating at least $(2\alpha_1+0.1\xi)k^{\prime}$ vertices, then by lemma~\ref{2-4}, there exists $T,U,S$ such that $C=T\bigcup U\bigcup S$,
$e_{F_1}(T,U)=0$, $d_{F_1[T]}(v)\leq 0.1\xi^6 k^{\prime}$ for every $v\in T$, and
\begin{align}\label{e2-1}
~~~~~~~~~~~~~~~~~~~~~~~~~~~~2|S|+|U|\leq (2\alpha_1+0.2\xi)k^{\prime}
\end{align}
for sufficient large $k^{\prime}$.

By~(\ref{e2-1}), $|T|\geq |C|-(2|S|+|U|)\geq(\alpha_2+9.8\xi)k^{\prime}$, $|S|\leq (\alpha_1+0.1\xi)k^{\prime}$,  and $|T|+|U|\geq
|C|-|S|\geq(2\alpha_1+\alpha_2+10\xi)k^{\prime}-(\alpha_1+0.1\xi)k^{\prime}=(\alpha_1+\alpha_2+9.9\xi)k^{\prime}$.

If $|T\bigcap V_1|\leq 3\xi k^{\prime}$ or $|T\bigcap V_2|\leq 3\xi k^{\prime}$, then, without loss of generality, assume that $|T\bigcap V_1|\leq 3\xi k^{\prime}$. Thus
$$|T\bigcap V_2|=|T|-|T\bigcap V_1|\geq(\alpha_2+6.8\xi)k^{\prime},$$
$$|V_1\backslash (S\bigcup T)|\geq |V_1|-|S|-|T\bigcap V_1|>(\alpha_2+4\xi)k^{\prime}.$$

\noindent Since $\xi^7 |V_1|^2\leq \xi^4|V_1\backslash (S\bigcup T)||T\bigcap V_2|$, then by Lemma~\ref{2-5}, $G[V_1\backslash (S\bigcup T), T\bigcap V_2]$ contains a blue connected
matching saturating at least $2(1-3\xi^4)\min\{|V_1\backslash (S\bigcup T)|$, $|T\bigcap V_2|\}> (2\alpha_2+7\xi)k^{\prime}$ vertices. Hence we may assume that $|T\bigcap
V_1|,|T\bigcap V_2|> 3\xi k^{\prime}$.

 All except at most $\xi^7 |V_1|^2+0.1\xi^6 k^{\prime}|T|<2\xi^6 (k^{\prime})^2$ pair $\{v,w\}$, $v\in T, w\in (V_1\bigcup V_2)\backslash S$ are blue edges of $G$. Note that
 $|V_1\backslash S|,|V_2\backslash S|\geq |V_1|-|S|> (\alpha_2+7\xi)k^{\prime}$. Thus

 $$2\xi^6 (k^{\prime})^2\leq \xi^4 \min\{|T\bigcap V_1||T\bigcap V_2|,|T\bigcap V_1||V_2\backslash S|,|T\bigcap V_2||V_1\backslash S|\}.$$

\noindent By lemma~\ref{2-5}, we can take $T_1\subseteq T\bigcap V_1,T_2\subseteq T\bigcap V_2, |T_1|\geq |T\bigcap V_1|- 2\xi^3 k^{\prime}, |T_2|\geq |T\bigcap V_2|-2\xi^3
k^{\prime}$, $U_1\subseteq V_1\backslash S, U_2\subseteq V_2\backslash S, |U_1|\geq |V_1\backslash S|-\xi^3 k^{\prime},|U_2|\geq |V_2\backslash S|-\xi^3 k^{\prime}$ such
that $G[T_1,T_2]$ contains at least a blue edge and $G[T_1,U_2],G[T_2,U_1]$ contain a blue connected matching saturating at least $2\min\{|T_1|,|U_2|\}$ and $2\min\{|T_2|,|U_1|\}$
vertices, respectively. Note that $|U_1|>|V_1\backslash S|-\xi^3 k^{\prime}>(\alpha_2+6\xi)k^{\prime}$,$|U_2|>|V_2\backslash S|-\xi^3 k^{\prime}>(\alpha_2+6\xi)k^{\prime}$. If
$|T_1|\geq |U_2|$ or $|T_2|\geq |U_1|$, then either $G[T_1,U_2]$ or $G[T_2,U_1]$ contains a blue connected matching saturating at least $(2\alpha_2+12\xi)k^{\prime}$ vertices.
Otherwise, $|T_1|< |U_2|$ and $|T_2|<|U_1|$. Also note that $|T_2|+|T_1|\geq |T|-4\xi^3 k^{\prime}>(\alpha_2+8\xi)k^{\prime}$. Thus $G[T_1\bigcup U_1,T_2\bigcup U_2]$ contains a blue
connected matching saturating at least $|T_2|+|T_1|>(\alpha_2+8\xi)k^{\prime}$ vertices.\q

\section{Multicolor bipartite Ramsey number of cycles}

Similar to that crucial part of the proof of Theorem~\ref{1-1} is to show  Lemma~\ref{2-3}, to complete the proof of Theorem~\ref{1-3}, it is crucial to show Lemma~\ref{3-1}. Next we
give the proof of Lemma~\ref{3-1}, and omit details of the proof of Theorem~\ref{1-3}.

\begin{lemma}\label{3-1}

For $r\geq 2$, $0<\alpha_1, \alpha_2\leq 1$, $[(j+2)!-1]\sum^{j+1}_{i=1} \alpha_i\leq \alpha_{j+2}\leq 1, j=1,2,\ldots,r-2$ and $0<\xi\leq
 \min\{((r+2)r!)^{-3}\alpha_1, ((r+2)r!)^{-3}\alpha_2, 10^{-3}\}$, there exist $N_0$ such that for each bipartite graph with bipartition $\{V_1,V_2\}$, $|V_1|=|V_2|=[\sum^r_{j=1}
\alpha_j+(r+2)r!\xi]k^{\prime}\geq N_0$ and
$|E(G)|\geq (1-\xi^{3r+5})|V_1|^2$, the following holds. For every $r$-edge coloring of $G$ there exists a color $i$, $i=1,2,\ldots,r$,
such that some component of the subgraph induced in $G$ by the edges of the $i$-th color contains a matching saturating at least $(2\alpha_i+0.1\xi)k^{\prime}$
vertices.
\end{lemma}

\p
 For $r\geq 2$, $0<\alpha_1, \alpha_2\leq 1$, $[(j+2)!-1]\sum^{j+1}_{i=1} \alpha_i\leq \alpha_{j+2}\leq 1, j=1,2,\ldots,r-2$ and $0<\xi\leq
 \min\{((r+2)r!)^{-3}\alpha_1, ((r+2)r!)^{-3}\alpha_2, 10^{-3}\}$. Let $G$ be a bipartite graph with bipartition $\{V_1,V_2\}$, $|V_1|=|V_2|=[\sum^r_{j=1}
 \alpha_j+(r+2)r!\xi]k^{\prime}$ and $|E(G)|\geq (1-\xi^{3r+5})|V_1|^2$. Let us consider an $r$-edge coloring of $G=\bigcup^r_{i=1}G_i$.

We proceed by induction on $r$. For $r=2$, by Lemma~\ref{2-3}, we are done. Now suppose that the lemma is true for $r-1$ and
prove it for $r$. Since $r!\xi k^{\prime}|V_1|>\xi^{3r+5}|V_1|^2$ and $|E(G)|\geq (1-\xi^{3r+5})|V_1|^2$, there exists a vertex $v\in V_2$ such that the number of the neighbors of
$v$ in $V_1$ is at least $[\sum^r_{j=1} \alpha_j+(r+1)r!\xi]k^{\prime}$. Thus there exists a monochromatic star $T$ with the center $v$ and at least $\frac{(\sum^r_{j=1}
\alpha_j+(r+1)r!\xi)k^{\prime}}{r}\geq [(r-1)!\sum^{r-1}_{j=1} \alpha_j+(r+1)(r-1)! \xi]k^{\prime}$ leaves.  Let $1\leq i\leq r$ and $H$ be a maximal component colored by the
$i$-th color with $C=V(H)$ and $|C\bigcap V_s|\geq [(r-1)!\sum^{r-1}_{j=1} \alpha_j+(r+1)(r-1)! \xi]k^{\prime}$ for some $s=1,2$. Without loss of generality, assume that
$$|C\bigcap
V_1|\geq [(r-1)!\sum^{r-1}_{j=1} \alpha_j+(r+1)(r-1)! \xi]k^{\prime}.$$
If $H$ contains a connected matching saturating at least $(2\alpha_i+0.1\xi)k^{\prime}$ vertices, then we are
done. Hence we may assume that $H$ does not contain a connected matching saturating at least $(2\alpha_i+0.1\xi)k^{\prime}$ vertices.

By Lemma~\ref{2-4}, there exists disjoint $T, U$ and $S$ such that $C=T\bigcup U\bigcup S$,
\begin{align}\label{e3-1}
~~~~~~~~~~~~~~~~~~~~~~~~~~~~2|S|+|U|\leq (2\alpha_i+0.2\xi)k^{\prime}
\end{align}
for sufficient large $k^{\prime}$, and in $H$, there are at most $\sqrt{[\sum^r_{j=1} \alpha_j+(r+2)r!\xi]k^{\prime}}<\xi^{3r+5} k^{\prime}$  neighbors in $T$ and no  neighbors in
$U$ for any $v\in T$. Clearly, $|S|\leq (\alpha_i+0.1\xi)k^{\prime}$.

\vskip3mm

\textbf{Case 1} $i=r$ and $|C\bigcap V_1|\geq [\sum^{r-1}_{j=1} a_j+(r+1)(r-1)! \xi]k^{\prime}$.

\vskip3mm

In this case, we can weaken the condition on $|C\bigcap V_1|$ to $|C\bigcap V_1|\geq [\sum^{r-1}_{j=1} a_j+(r+1)(r-1)! \xi]k^{\prime}$. Note that
$\xi^{3r+5}|V_1|^2<\xi^{3r+2}[(\sum^{r-1}_{j=1} \alpha_j+(r+1)(r-1)!
\xi)k^{\prime}]^2$. If $|V_2\backslash C|\geq (\sum^{r-1}_{j=1} \alpha_j+(r+1)(r-1)! \xi)k^{\prime}$, then by the induction hypothesis,  $G[C\bigcap V_1, V_2\backslash C]$ contains
a matching saturating at least $(2\alpha_j+0.1\xi)k^{\prime}$
vertices colored by the $j$-th color for some $1\leq j\leq r-1$. Hence we may assume that $|V_1\bigcap C|,|V_2\bigcap C|\geq [\alpha_r+(r^2+r-1)(r-1)! \xi]k^{\prime}$.

By~(\ref{e3-1}), either $| V_1\bigcap (U\bigcup S)|\leq (a_r+0.1\xi)k^{\prime}$ or $|V_2\bigcap (U\bigcup S)| \leq (a_r+0.1\xi)k^{\prime}$.
Without loss of generality, assume that $|V_1\bigcap (U\bigcup S)| \leq (a_r+0.1\xi)k^{\prime}$. Thus
$$|V_1\backslash(U\bigcup S)|=|V_1|-|V_1\bigcap (U\bigcup S)|> (\sum^{r-1}_{j=1} a_j+(r+1)(r-1)! \xi)k^{\prime}.$$

Note that $\xi^{3r+5} (\sum^{r-1}_{j=1} a_j+(r+1)(r-1)! \xi)(k^{\prime})^2+\xi^{3r+5}|V_1|^2< \xi^{3r+2}[(\sum^{r-1}_{j=1} a_j+(r+1)(r-1)! \xi)k^{\prime}]^2$. If
$|V_2\bigcap(U\bigcup T)|\geq (\sum^{r-1}_{j=1} a_j+(r+1)(r-1)! \xi)k^{\prime}$, then  by the induction hypothesis,  $G[V_2\bigcap(U\bigcup T), V_1\backslash(U\bigcup S)]$ contains a
matching saturating at least  $(2\alpha_j+0.1\xi)k^{\prime}$ vertices colored by the $j$-th color for some $1\leq j\leq r-1$. Since $|S|\leq (a_r+0.1\xi)k^{\prime}$, $|V_2\backslash
S|>(\sum^{r-1}_{j=1} a_j+(r+1)(r-1)! \xi)k^{\prime}$. If $|T\bigcap V_1|\geq (\sum^{r-1}_{j=1} a_j+(r+1)(r-1)! \xi)k^{\prime}$, then, similarly, $G[V_2\backslash S, T\bigcap V_1]$
contains a matching saturating at least $(2\alpha_j+0.1\xi)k^{\prime}$ vertices colored by the $j$-th color for some $1\leq j\leq r-1$.

Hence we may assume that $|T\bigcap V_1|,|V_2\bigcap(U\bigcup T)|<(\sum^{r-1}_{j=1} a_j+(r+1)(r-1)! \xi)k^{\prime}$. Clearly, $|V_2\bigcap S|=|C\bigcap V_2|-|V_2\bigcap(U\bigcup
T)|>[a_r-\sum^{r-1}_{j=1} a_j+(r^2-2)(r-1)! \xi)]k^{\prime}$ and $|V_1\bigcap(U\bigcup S)|=|C\bigcap V_1|-|T\bigcap V_1|>[a_r-\sum^{r-1}_{j=1} a_j+(r^2-2)(r-1)! \xi]k^{\prime}$.
Since $|S|\leq (a_r+0.1\xi)k^{\prime}$, then $|V_1\bigcap S|<\{\sum^{r-1}_{j=1} a_j+[0.1-(r^2-2)(r-1)! ] \xi\}k^{\prime}$. Thus
\begin{align*}
|V_1\bigcap (T\bigcup U)|&=|V_1\bigcap C|-|V_1\bigcap S|\\
&>\{(r-1)!\sum^{r-1}_{j=1} a_j+[(2r^2+r+1)(r-1)! -0.1] \xi\}k^{\prime}\\
&>(\sum^{r-1}_{j=1} a_j+(r+1)(r-1)! \xi)k^{\prime}
\end{align*}
\noindent and
\begin{align*}
|V_2\bigcap (U\bigcup S)|&<|U|+2|S|-(|V_2\bigcap S|+|V_1\bigcap (U\bigcup S)|)\\
&<(2a_r+0.2\xi)k^{\prime}-2[a_r-\sum^{r-1}_{j=1} a_j+(r^2-2)(r-1)! \xi)]k^{\prime}\\
&<2\sum^{r-1}_{j=1} a_j.
\end{align*}

Thus $|V_2\backslash(U\bigcup S)|=|V_2|-|V_2\bigcap (U\bigcup S)|>(\sum^{r-1}_{j=1} a_j+(r+1)(r-1)! \xi)k^{\prime}$. By the induction hypothesis, $G[V_2\backslash (U\bigcup S),
V_1\bigcap (T\bigcup U)]$ contains a matching saturating at least  $(2\alpha_j+0.1\xi)k^{\prime}$ vertices colored by the $j$-th color for some $1\leq j\leq r-1$.

\vskip3mm

\textbf{Case 2} $1\leq i\leq r-1$.

\vskip3mm

Recall that $|C\bigcap V_1|\geq [(r-1)!\sum^{r-1}_{j=1} \alpha_j+(r+1)(r-1)! \xi]k^{\prime}$. If $|V_2\backslash C|\geq [(r-1)!\sum^{r-1}_{j=1} \alpha_j+(r+1)(r-1)! \xi]k^{\prime}$,
then, by the
induction hypothesis, $G[C\bigcap V_1, V_2\backslash C]$ contains a connected matching saturating at least $(2\alpha_j+0.1\xi)k^{\prime}$ vertices colored by
the $j$-th color for some $1\leq j\leq r-1, j\neq i$ or $\{2[(r-1)!-1]\sum^{r-1}_{j=1} \alpha_j+2\alpha_i+0.1\xi\}k^{\prime}$ vertices colored by the $r$-th color. Since
$\alpha_1,\alpha_2\geq ((r+2)r!)^{3}\xi$ and $[(j+2)!-1]\sum^{j+1}_{i=1} \alpha_i\leq \alpha_{j+2},j=1,2,\ldots,r-2$, $\{[(r-1)!-1]\sum^{r-1}_{j=1}
\alpha_j+\alpha_i\}k^{\prime}>[\sum^{r-1}_{j=1} \alpha_j+(r+1)(r-1)! \xi]k^{\prime}$ . Combine with Case 1, the proof is complete. Hence we can assume that $|V_2\backslash C|<
[(r-1)!\sum^{r-1}_{j=1} \alpha_j+(r+1)(r-1)! \xi]k^{\prime}$.

Clearly,  $|V_2\bigcap C|> [(r-1)!(r-1)\sum^{r-1}_{j=1} \alpha_j+(r^2+r-1)(r-1)! \xi]k^{\prime}$. Since $|V_2\bigcap S|<(\alpha_i+0.1\xi)k^{\prime}$,
$|V_2\bigcap(C\backslash S)|>\{(r-1)!+1]\sum^{r-1}_{j=1} \alpha_j+2(r+1)(r-1)! \xi\}k^{\prime}$. By~(\ref{e3-1}), $|V_1\backslash(U\bigcup S)|\geq [\sum^r_{j=1}
\alpha_j+(r+2)r!\xi]k^{\prime}-(2\alpha_i+0.2\xi)k^{\prime}>\{(r-1)!+1]\sum^{r-1}_{j=1} \alpha_j+2(r+1)(r-1)! \xi\}k^{\prime}$.

Let $\beta_1=\alpha_1,\ldots, \beta_{i-1}=\alpha_{i-1}, \beta_{i+1}=\alpha_{i+1},\ldots,\beta_{r-1}=\alpha_{r-1}, \beta_{r}=[(r-1)!\sum^{r-1}_{j=1} \alpha_j+\alpha_i+(r+1)(r-1)!
\xi]k^{\prime}$. By the
induction hypothesis, $G[V_2\bigcap(C\backslash S), V_1\backslash (U\bigcup S)]$ contains a connected matching saturating at least $(2\beta_j+0.1\xi)k^{\prime}$ vertices colored by
the $j$-th color for some $1\leq j\leq r, j\neq i$. If $j\in [r]\backslash\{i,r\}$, then the proof is complete. If $j=r$, then $G[V_2\bigcap(C\backslash S), V_1\backslash (U\bigcup
S)]$ contains a connected matching saturating at least $(2\beta_r+0.1\xi)k^{\prime}$ vertices colored by the $r$-th color. By the discussion of Case 1, we are done.

\q

\section{2-colored bipartite graphs with ``large" minimum degree}

Denote by $G_{R}[H]$ and $G_{B}[H]$the induced subgraph of red edges and blue edges of $G$ with respect to $H\subseteq V(G)$, respectively. Denote by  $N_{R}(v)$ (or $N_{B}(v)$) the
set of red (or blue) neighbors of a vertex$v$ in $G$, respectively.

\subsection{Degree form of the Regularity Lemma}

In the proof of Theorem~\ref{1-5}, We shall use the degree form of the Szemer\'edi Regularity Lemma for an edge-colored bipartite graphs (see \cite{JK96},
Theorem 1.10).

\begin{lemma}\emph{(}Degree form of Red-Blue Regularity Lemma\emph{)}\label{4-1}
For every $\epsilon>0$ and positive integer $k_0\geq \frac{1}{\epsilon}$, there exists $M=M(\epsilon,k_0)$ such that for any $d\in[0,1]$ and any $2$-edge coloring of a bipartite
graph with
$|X|=|Y|=N$, there is a partition of $X$ into clusters $X_0,X_1,\ldots,X_k$, a partition of $Y$ into clusters $Y_0,Y_1,\ldots,Y_k$, and a subgraph $G^{\prime}\subseteq G$ with the
following properties:
\begin{enumerate}[(1)]
\item $k_0\leq k\leq M$;
\item $|X_0|=|Y_0|\leq \epsilon N$;
\item $|X_i|=|Y_i|=m\leq \epsilon N$ for all $i\geq 1$;
\item $d_{G^{\prime}} (v)>d_{G} (v)-(2d+\epsilon)N$ for all $v\notin X_0 \bigcup Y_0$;
\item for all $1\leq i,j\leq k$, the pair $(X_i,Y_j)$ is $\epsilon$-regular for $G_R[V(G^{\prime})]$ with a density either $0$ or greater than $d$ and  $\epsilon$-regular for
    $G_B[V(G^{\prime})]$ with a density either $0$ or greater than $d$, where $E(G^{\prime})=E(G_R[V(G^{\prime})])\bigcup E(G_B[V(G^{\prime})])$ is the induced $2$-edge colouring
    of $G^{\prime}$.
\end{enumerate}
\end{lemma}

Having applied the above form of the Regularity Lemma to a $2$-edge coloured bipartite graph, we need the following definition.

\begin{defi}\emph{(}$(\epsilon, d)$-reduced Graph\emph{)}\label{4-2}
Given a bipartite graph $G=(X,Y;E)$ and a partition $\{X_i: 1\leq i \leq k\}$ of $X$ and $\{Y_i: 1\leq i \leq k\}$ of $Y$ from Lemma~\ref{4-1}, we define the $(\epsilon, d)$-reduced
$2$-coloured bipartite graph $H$ on vertex set $\{v^{(s)}_{i}:1\leq i\leq k, s=1,2\}$ as follows: For
$1\leq i,j\leq k$,
\begin{enumerate}[(1)]
\item let $(v^{(1)}_{i},v^{(2)}_{j})$ be an edge of $H$ when $G_{B}[X_{i},Y_{j}]$ has density at least $d$;
\item let $(v^{(1)}_{i},v^{(2)}_{j})$ be an edge of $H$ when it is not a blue edge and $G_{R}[X_{i},Y_{j}]$ has density at least $d$.
\end {enumerate}
\end{defi}

\subsection{Proof of Theorem~\ref{1-5}}

\vskip3mm

We may assume that $0<\xi\leq 10^{-3}$. Let $G$ be a bipartite graph with bipartition $\{V_1,V_2\}$, $|V_1|=|V_2|=N=(2+8\xi)n$, and
$\delta(G)>(\frac{7}{8}+9\xi)N$. Consider a $2$-edge coloring of $G$, i.e., $G_1, G_2$. Choose $\epsilon=\xi^3$ and $d=\xi$. Apply
Lemma~\ref{4-1} to $G$, with parameters $d$ and $\epsilon$. Let $V^{(s)}_0,V^{(s)}_1,\ldots,V^{(s)}_k$ be the partition of bipartition $V^{(s)}$ for $s=1,2$
and $G^{\prime}$ be the subgraph of $G$ guaranteed by Lemma~\ref{4-1}. Finally,
let $H$ be the $(\epsilon, d)$-reduced graph deduced from $G^{\prime}$, with $2$-edge colouring $H=H_R\bigcup H_B$.

We first observe that
$$\delta(H)\geq (\frac{7}{8}+\xi)k.$$
Indeed, by (4), we have $\delta(G^{\prime})\geq (\frac{7}{8}+8\xi-2d)N$. Suppose that $\delta(H)<(\frac{7}{8}+\xi)k$. Then there exists some $i\geq 1$ and some $s\in\{1,2\}$
such that $d_{H}(v^{(s)}_i)<(\frac{7}{8}+\xi)k$. For a vertex $v\in V^{(s)}_i$, its neighbours in $G^{\prime}$ are only in $V^{(3-s)}_0$, or in $V^{(3-s)}_j$
for
those $j$ such that $v^{(s)}_iv^{(3-s)}_j$ is an edge of $H$. Hence
$$\delta(G^{\prime})<(\frac{7}{8}+\xi)km+|V^{(3-s)}_0|<(\frac{7}{8}+2\xi)N,$$
contradicting to $\delta(G^{\prime})\geq (\frac{7}{8}+8\xi-2d)N$. Then we will prove the following crucial lemma whose proof will be given in subsection 4.3.

\begin{lemma}\label{4-3}
For $0<\xi\leq 10^{-3}$, there exist
$N_0$, such that for each bipartite graph with bipartition $\{V_1,V_2\}$, $|V_1|=|V_2|=(2+8\xi)k^{\prime}\geq N_0$ and
$\delta(G)\geq(\frac{7}{8}+\xi)(2+8\xi)k^{\prime}$, the following holds. For any $2$-edge coloring of $G$, there exists a monochromatic component  containing a matching saturating
at least $(2+0.1\xi)k^{\prime}$
vertices.
\end{lemma}

Apply Lemma~\ref{4-3} with $k=(2+8\xi)k^{\prime}$ to the $2$-edge coloring of $H=H_R\bigcup H_B$. Thus there exists some $f\in\{R,B\}$, such that one of the components of $H_f$
contains a matching $M^*$ saturating at least $2a=(2+0.1\xi)k^{\prime}$ vertices. Let $F^*$ be a minimal
connected tree of $H_f$ containing $M^*$. Consider a closed walk $W=v^{(s)}_{i_1}v^{(3-s)}_{i_2}\cdots v^{(3-s)}_{i_t}v^{(s)}_{i_1}$ which contains all edges of $M^*$. Since $F^*$
is a tree, so $W$ must be of even length.

Applying Fact~\ref{2-0} repeatedly, we can show the existence of  $v^{(q)}_{i_p}\in V^{(q)}_{i_p}$, $1\leq p\leq t+1, q=s$ or $3-s$, $i_{t+1}=i_1$, satisfying:

\begin{enumerate}[(1)]
  \item $v^{(q)}_{i_p}$ has at least $\frac{4}{5}\epsilon|V^{(q)}_{i_p}|$ neighbors in both $V^{(3-q)}_{i_{p-1}}$ and $V^{(3-q)}_{i_{p+1}}$ for each $p$;
  \item If $v^{(q)}_{i_p}v^{(3-q)}_{i_{p+1}}$ is not an edge in $M^{*}$(view an edge in the closed walk as in $M^{*}$ only when it is in $M^{*}$ and the closed walk first passes
      it), then $v^{(q)}_{i_p}v^{(3-q)}_{i_{p+1}}$ is an edge in $H_f$.
\end{enumerate}

Let $m=\lfloor\frac{N}{k}\rfloor$.
Applying Lemma~\ref{2-1} with $\beta=4\xi$, for every $l$, $1\leq l\leq (1-\frac{5\xi^2}{4}) m$, each edge $v^{(q)}_{i_p}v^{(3-q)}_{i_{p+1}}$ in $M^*$ can be enlarged to a path of
length $2l+1$ connecting $u^{(q)}_{i_p}\in V^{(q)}_{i_p}$
and $u^{(3-q)}_{i_{p+1}}\in V^{(3-q)}_{i_{p+1}}$. Therefore,
there exists a cycle in $G_f$ of each even length $l^{\prime}= \sum^{a}_{j=1}2l_j+t$, where $1\leq l_j\leq (1-\frac{5\xi^2}{4}) m$ for $1\leq j\leq a$,  and
\begin{align*}
\sum^{a}_{j=1}2(1-\frac{5\xi^2}{4}) m+t&\geq 2a(1-\frac{5\xi^2}{4}) \frac{N(1-\epsilon)}{k}+t\\
&\geq(2+0.1\xi)k^{\prime}(1-2\xi^2)(1-\xi^3)\frac{N}{k}+t\\
&>(2+0.05\xi)n+t.
\end{align*}
Therefore, $G_f$ contains a cycle of length $2n$ for some $f\in \{R,B\}$.

\subsection{Proof of Lemma~\ref{4-3}}

\begin{lemma}[\cite{LD16}]\label{4-4}
Let $n$ be even and $k\geq 2$, and let $G$ be a $k$-partite graph on $n$ vertices with the vertex set partitioned as $\{X_1,X_2,\ldots,X_k\}$. Suppose that
$|X_i|\leq \frac{n}{2}$ for all $i\in [k]$. If $d(x)>\frac{3}{4}n-|X_i|$ for all $i\in [k]$ and for all $x\in X_i$, then $G$ is connected and contains a perfect matching.
\end{lemma}

The following corollary is immediately  from  Lemma~\ref{4-4}.

\begin{corollary}\label{4-5}
Let $G$ be a bipartite graph with bipartition $\{X,Y\}$, $|X|=|Y|=n$. If $d(x)>\frac{1}{2}n$ for $x\in X\bigcup Y$, then $G$ is connected and contains a perfect matching.
\end{corollary}

We will show the following result.

\begin{lemma}\label{4-6}
Let $0<\xi\leq 10^{-3}$. Let $G$ be a bipartite graph with bipartition $\{V_1,V_2\}$, $|V_1|=|V_2|=(2+8\xi)k^{\prime}$ and $\delta(G)>(\frac{3}{4}+\xi)(2+8\xi)k^{\prime}$.
Then any $2$-edge coloring of $G$ leads to a monochromatic component with at least $(2+6\xi)k^{\prime}$ vertices.
\end{lemma}

\p Let us consider a $2$-edge coloring of $G$ and let $F_1$ be a monochromatic component of $G$ with the largest number of vertices (say, red, and another color is blue). Let
$C=V(F_1)$.
If $|C\bigcap V_1|\geq (1+3\xi)k^{\prime}$ or $|C\bigcap V_2|\geq (1+3\xi)k^{\prime}$, then, without loss of generality, assume that $|C\bigcap V_1|\geq (1+3\xi)k^{\prime}$ and
$|C\bigcap V_2|<(1+3\xi)k^{\prime}$. Otherwise, $|C|\geq (2+6\xi)k^{\prime}$ and we are done. Thus $|V_2\backslash C|>(1+3\xi)k^{\prime}$. Take $H_1\subseteq C\bigcap V_1
,H_2\subseteq V_2\backslash C$ such that $H_1=H_2=(1+3\xi)k^{\prime}$. Since $\delta(G)>(\frac{3}{4}+\xi)(2+8\xi)k^{\prime}$, we have
$$\delta(G_{B}[H_1,H_2])\geq (\frac{3}{4}+\xi)(2+8\xi)k^{\prime}-(1+5\xi)k^{\prime} >\frac{1}{2}\times (1+3\xi)k^{\prime}.$$
By  Corollary~\ref{4-5}, $G_{B}([H_1,H_2])$ is connected and has $(2+6\xi)k^{\prime}$ vertices. Hence we may assume that $|C\bigcap V_1|,|C\bigcap V_2|<(1+3\xi)k^{\prime}$. For any
$u\in C\bigcap V_1$,
\begin{align*}
|N_{B}(u)\bigcap (V_2\backslash C)|
&\ge (\frac{3}{4}+\xi)|V_2|-|C\bigcap V_2|\\
&=(\frac{3}{4}+\xi)(|V_2\backslash C|+|C\bigcap V_2|)-|C\bigcap V_2|\\
&>\frac{3}{4}|V_2\backslash C|-\frac{1}{4}|C\bigcap V_2|\\
&>\frac{1}{2}|V_2\backslash C|.
\end{align*}
So $G_{B}[C\bigcap V_1, V_2\backslash C]$ is blue and connected, and this is a large monochromatic component than $F_1$, a contradiction.
\q

Now we will give the proof of Lemma~\ref{4-3}.

\vskip3mm

\noindent\textbf{Proof of Lemma~\ref{4-3}.} Let $G$ be a bipartite graph with bipartition $\{V_1,V_2\}$, $|V_1|=|V_2|=(2+8\xi)k^{\prime}\geq N_0$ and
$\delta(G)\geq(\frac{7}{8}+\xi)(2+8\xi)k^{\prime}$. Let us consider a $2$-edge coloring of $G$, say, red and blue. Let $F_1$ denote a monochromatic component of $G$ with the
largest number
of vertices and $C=V(F_1)$. Without loss of generality, assume that the edges of $F_1$ are colored by red and $F_1$ does not contain a matching saturating at least
$(2+0.1\xi)k^{\prime}$ vertices. By Lemma~\ref{4-6}, we have
$$|C|\geq (2+6\xi)k^{\prime}.$$

By Lemma~\ref{2-4}, there exists disjoint $T, U$ and $S$ such that $C=T\bigcup U\bigcup S$,

\begin{align}\label{e4-1}
~~~~~~~~~~~~~~~~~~~~~~~~~~~~2|S|+|U|\leq (2+0.2\xi)k^{\prime}
\end{align}
for sufficient large $k^{\prime}$, and for any $v\in T$, there are at most $\sqrt{(2+8\xi)k^{\prime}}<\xi^2 k^{\prime}$ red neighbors in $T$ and no red neighbors in $U$ . Clearly,
$|S|\leq (1+0.1\xi)k^{\prime}$.

\begin{claim}\label{4-9}
If $|(T\bigcup U)\bigcap V_s|$, $|V_{3-s}\backslash (U\bigcup S)|\geq(1+4\xi)k^{\prime}$ for some $s\in\{1,2\}$, then $G_{B}[(T\bigcup U) \bigcap V_s ,V_{3-s}\backslash (U\bigcup
S)]$
contains a connected matching saturating at least $(2+8\xi)k^{\prime}$ vertices.
\end{claim}

\p Take $H_1\subseteq (T\bigcup U)\bigcap V_s, H_2\subseteq V_{3-s}\backslash (U\bigcup S)$, $|H_1|=|H_2|=(1+4\xi)k^{\prime}$. Then
$$\delta(G_{B}[H_1,H_2])\geq (\frac{7}{8}+\xi)(2+8\xi)k^{\prime}-(1+4\xi)k^{\prime}-\xi^2 k^{\prime}>\frac{1}{2}(1+4\xi)k^{\prime}.$$
By Corollary~\ref{4-5}, $G_{B}[H_1,H_2]$ contains a connected matching saturating at least $(2+8\xi)k^{\prime}$ vertices. \q

\begin{claim}\label{4-7} If $ |C\bigcap V_s|, |V_{3-s}\backslash C|\geq (1+4\xi)k^{\prime}$ for some $s\in\{1,2\}$, then
$G_{B}[C\bigcap V_s, V_{3-s}\backslash C]$ contains a connected matching saturating at least $(2+8\xi)k^{\prime}$ vertices.
\end{claim}
\p Take any $H_1\subseteq C\bigcap V_s, H_2\subseteq  V_{3-s}\backslash C, |H_1|=|H_2|=(1+4\xi)k^{\prime}$. Then
$$\delta (G_{B}[H_1,H_2])\geq (\frac{7}{8}+\xi)(2+8\xi)k^{\prime}-(1+4\xi)k^{\prime}> \frac{1}{2}(1+4\xi)k^{\prime}.$$
By Corollary~\ref{4-5}, $G_{B}[H_1,H_2]$ contains a connected matching saturating at least $(2+8\xi)k^{\prime}$ vertices.
\q

\begin{claim}\label{4-8} If $|V_s\backslash C|\geq \frac{1}{2}(1+4\xi)k^{\prime},|C\bigcap V_s|\geq (1+4\xi)k^{\prime}$ for each $s=1,2$, then $G_{B}[V_1 ,V_2]$ contains a
connected matching saturating at least $(2+8\xi)k^{\prime}$ vertices.
\end{claim}
\p First we show that $e(G_{B}[C\bigcap V_1,C\bigcap V_2])\geq 1$. By contradiction, suppose that $e(G_{B}[C\bigcap V_1,C\bigcap V_2])=0$. Take any $H_1\subseteq C\bigcap V_1,
H_2\subseteq C\bigcap V_2$ such that $|H_1|=|H_2|=(1+4\xi)k^{\prime}$. Thus
$$\delta (G_{R}[H_1,H_2])\geq (\frac{7}{8}+\xi)(2+8\xi)k^{\prime}-(1+4\xi)k^{\prime}> \frac{1}{2}(1+4\xi)k^{\prime}.$$
By Corollary~\ref{4-5}, $G_{R}[H_1,H_2]$ contains a connected matching saturating at least $(2+8\xi)n$ vertices, a contradiction to the hypothesis of $F_1$. Thus there exists an
edge $uv\in E(G_{B}[C\bigcap V_1,C\bigcap V_2])$, where $u\in C\bigcap V_1, v\in C\bigcap V_2$.

Now take any $H_1\subseteq V_1\backslash C , H_2\subseteq V_2\backslash C,D_1\subseteq  C\bigcap V_1, D_2\subseteq  C\bigcap V_2$ such that $u\in D_1$,  $v\in D_2$ and
$|H_1|=|H_2|=|D_1|=|D_2|=\frac{1}{2}(1+4\xi)k^{\prime}$. Then
$$\delta (G_{B}[H_1,D_2]),\delta (G_{B}[H_2,D_1])\geq (\frac{7}{8}+\xi)(2+8\xi)k^{\prime}-\frac{3}{2}(1+4\xi)k^{\prime}> \frac{1}{2}\times[\frac{1}{2}(1+4\xi)k^{\prime}].$$
By Corollary~\ref{4-5}, $G_{B}[H_1,D_2]$ and $G_{B}[H_2,D_1]$ are connected and contain a perfect matching. Note that $G_{B}[D_1,D_2]$ contains an edge $uv$. Thus $G_{B}[H_1\bigcup
D_1,H_2\bigcup D_2]$ contains a connected matching saturating at least $(2+8\xi)k^{\prime}$ vertices.\q

 If $|V_1\bigcap C|< (1+4\xi)k^{\prime}$ or $|V_2\bigcap C|<(1+4\xi)k^{\prime}$, then, without loss of generality, assume that $|V_1\bigcap C|<
 (1+4\xi)k^{\prime}$. Clearly, $|V_1\backslash C|=|V_1|-|V_1\bigcap C|> (1+4\xi)k^{\prime}$. Note that $|C|\geq (2+6\xi)k^{\prime}$. Thus $|V_2\bigcap
 C|=|C|-|V_1\bigcap C|>(1+2\xi)k^{\prime}.$
By Claim~\ref{4-7}, $G_{B}[V_1\backslash C,V_2\bigcap C]$ contains a connected matching saturating at least $(2+4\xi)k^{\prime}$ vertices. Hence we may assume that
 $|V_s\bigcap C|\geq(1+4\xi)k^{\prime}$,
$|V_{3-s}\backslash C|<(1+4\xi)k^{\prime}$ for each $s=1,2$. 
 By Claim~\ref{4-8}, we may assume
that either $|V_1\backslash C|<\frac{1}{2}(1+4\xi)k^{\prime}$ or $|V_2\backslash C|< \frac{1}{2}(1+4\xi)k^{\prime}$.
Without loss of generality, assume that $|V_1\backslash C|<\frac{1}{2}(1+4\xi)k^{\prime}$. Then $|V_1\bigcap C|> \frac{3}{2}(1+4\xi)k^{\prime}$.

To summary, from now on,  we may assume that  $|V_2\bigcap C|\geq(1+4\xi)k^{\prime}$,
$|V_2\backslash C|<(1+4\xi)k^{\prime}$,  $|V_1\backslash C|<\frac{1}{2}(1+4\xi)k^{\prime}$, and  $|V_1\bigcap C|> \frac{3}{2}(1+4\xi)k^{\prime}$.

\vskip3mm

\textbf{Case 1} $|C|\geq(3+10\xi)k^{\prime}$.

\vskip3mm

By~(\ref{e4-1}), $|S|\leq (1+0.1\xi)k^{\prime}$, $|T|+|U|\geq |C|-|S|\geq(2+9.9\xi)k^{\prime}$.

\vskip3mm

\textbf{Subcase 1.1} $|(T\bigcup U)\bigcap V_1|, |(T\bigcup U)\bigcap V_2|\geq(1+4\xi)k^{\prime}$.

\vskip3mm

By~(\ref{e4-1}), either $|U\bigcap V_1|+|S\bigcap V_1|\leq (1+0.1\xi)k^{\prime}$ or $|U\bigcap V_2|+|S\bigcap V_2|\leq (1+0.1\xi)k^{\prime}$. Without loss of generality, we may
assume
that $|U\bigcap V_1|+|S\bigcap V_1|\leq (1+0.1\xi)k^{\prime}$. Thus
$$|V_1\backslash (U\bigcup S)|=|V_1|-|(U\bigcup S)\bigcap V_1|\geq (1+7.9\xi)k^{\prime}.$$

By Claim~\ref{4-9}, $G_{B}[V_1\backslash (U\bigcup S),(T\bigcup U)\bigcap V_2]$ contains a connected matching saturating at least $(2+8\xi)k^{\prime}$ vertices.

\vskip3mm

\textbf{Subcase 1.2} $|(T\bigcup U)\bigcap V_1|<(1+4\xi)k^{\prime}$ or  $|(T\bigcup U)\bigcap V_2|<(1+4\xi)k^{\prime}$.

\vskip3mm

Without loss of generality, assume that $|(T\bigcup U)\bigcap V_2|<(1+4\xi)k^{\prime}$. Recall that $|T|+|U|\geq(2+9.9\xi)k^{\prime}$.
Thus
$$|(T\bigcup U)\bigcap V_1|=|T|+|U|-|(T\bigcup U)\bigcap V_2|> (1+5.9\xi)k^{\prime}.$$
If $|(U\bigcup S)\bigcap V_2|\leq (1+4\xi)k^{\prime}$, then
$$|V_2\backslash (U\bigcup S)|=|T\bigcap V_2|+|V_2\backslash C|=|V_2|-|(U\bigcup S)\bigcap V_2|\geq (1+4\xi)k^{\prime}.$$
By Claim~\ref{4-9}, $G_{B}[(T\bigcup U)\bigcap V_1,V_2\backslash (U\bigcup S)]$ contains a connected matching saturating at least $(2+8\xi)k^{\prime}$ vertices. Hence we may assume
that $|(U\bigcup S)\bigcap V_2|>(1+4\xi)k^{\prime}$. Since $|(T\bigcup U)\bigcap V_2|<(1+4\xi)k^{\prime}$, $|S\bigcap V_2|>|T\bigcap V_2|$, then
$$|T\bigcap V_1|=|C|-|U\bigcup S \bigcup (T\bigcap V_2)|\geq |C| -(|U|+2|S|)\geq (1+9.8\xi)k^{\prime}.$$
Note that $|V_2\backslash S|\geq |V_2|-|S|\geq (1+7.9\xi)k^{\prime}$. Similar to Claim~\ref{4-9}, $G_{B}[T\bigcap V_1,V_2\backslash S]$ contains a connected matching saturating at
least
$(2+8\xi)k^{\prime}$ vertices.
\q

\vskip3mm

\textbf{Case 2} $|C|<(3+10\xi)k^{\prime}$.

\vskip3mm

Recall  that $|V_1\backslash
C|+|V_2\backslash C|\geq (1+6\xi)k^{\prime}$,  $|V_2\backslash C|<(1+4\xi)k^{\prime}$  and $|V_1\backslash C|<{1 \over 2}(1+4\xi)k^{\prime}$. Thus
$$\frac{1}{2}(1+4\xi)k^{\prime}<|V_2\backslash C|< (1+4\xi)k^{\prime}, |V_2\bigcap C|>
(1+4\xi)k^{\prime},$$
$$|C|= 2|V_1|-(|V_1\backslash C|+|V_2\backslash C|)>(4+16\xi)k^{\prime}-\frac{3}{2}(1+4\xi)k^{\prime}=\frac{5}{2}(1+4\xi)k^{\prime}.$$

\begin{claim}\label{4-10} Under the conditions of Case 2,
if $|T\bigcap V_1|\geq \frac{1}{2}(1+4\xi)k^{\prime}$ or $|T\bigcap V_2|\geq \frac{1}{2}(1+4\xi)k^{\prime}$, then $G_{B}[C\bigcap V_1,V_2\backslash S]$ contains a connected
matching saturating $(2+8\xi)k^{\prime}$ vertices.
\end{claim}

\p  Since $|S|\leq (1+0.1\xi)k^{\prime}$, $|V_2\backslash S|> (1+4\xi)k^{\prime}$. Recall that $ |C\bigcap V_1|>\frac{3}{2}(1+4\xi)k^{\prime}$, $|V_2\backslash C|>
\frac{1}{2}(1+4\xi)k^{\prime}$. Then $|(T\bigcup U)\bigcap V_1|\geq {1 \over 2}(1+4\xi)k^{\prime}$. If $|T\bigcap V_1|\geq \frac{1}{2}(1+4\xi)k^{\prime}$ or $|T\bigcap V_2|\geq
\frac{1}{2}(1+4\xi)k^{\prime}$, then take $H_1\subseteq T\bigcap V_1$, $D_1\subseteq (V_1\bigcap C)\backslash H_1$, $D_2\subseteq V_2\backslash C$, $H_2\subseteq V_2\backslash
(S\bigcup D_2)$  or $H_1\in (T\bigcup U)\bigcap V_1$, $H_2 \in T\bigcap V_2 $, $D_1\subseteq (V_1\bigcap C)\backslash H_1,  D_2\subseteq V_2\backslash C$ such that
$|H_1|=|H_2|=|D_1|=|D_2|=\frac{1}{2}(1+4\xi)k^{\prime}$.  Thus
$$\delta(G_{B}[H_1,H_2])\geq (\frac{7}{8}+\xi)(2+8\xi)k^{\prime}-\frac{3}{2}(1+4\xi)k^{\prime}-\xi^2 k^{\prime}>\frac{1}{2}\times[\frac{1}{2}(1+4\xi)k^{\prime}].$$
$$\delta(G_{B}[D_1,D_2])\geq(\frac{7}{8}+\xi)(2+8\xi)k^{\prime}-\frac{3}{2}(1+4\xi)k^{\prime}>\frac{1}{2}\times[\frac{1}{2}(1+4\xi)k^{\prime}].$$
By Corollary~\ref{4-5}, $G_{B}[H_1,H_2]$ and $G_{B}[D_1,D_2]$ are connected and contain a perfect matching saturating $(1+4\xi)k^{\prime}$ vertices. Note that
$e(G_{B}[H_1,D_2])\neq 0$. Suppose not, then for all $v\in H_1$ or $D_2$, $d_{G}(v)\leq (2+8\xi)k^{\prime}-\frac{1}{2}(1+4\xi)k^{\prime}=\frac{3}{2}(1+4\xi)k^{\prime}$,
contradicting to
$\delta(G)\geq(\frac{7}{8}+\xi)(2+8\xi)k^{\prime}$. Thus $G_{B}[H_1\bigcup D_1,H_2\bigcup D_2]$ contain a connected matching saturating at least $(2+8\xi)k^{\prime}$ vertices.
\q

If $|(T\bigcup U)\bigcap V_1|\geq(1+4\xi)k^{\prime}$, then by Claim~\ref{4-9},
 we may assume that
$|V_2\backslash (U\bigcup S)|<(1+4\xi)k^{\prime}$. Then $|(U\bigcup S)\bigcap V_2|> (1+4\xi)k^{\prime}$. By~(\ref{e4-1}), $|(U\bigcup S)\bigcap V_1|<(1+4\xi)k^{\prime}$ and hence
$|V_1\backslash (U\bigcup S)|> (1+4\xi)k^{\prime}$. If
$|(U\bigcup T)\bigcap V_2|\geq (1+4\xi)k^{\prime}$, then, by Claim~\ref{4-9}, $G_{B}[V_1\backslash (U\bigcup S),(U\bigcup T)\bigcap V_2]$ contains a connected matching saturating
at least $(2+8\xi)k^{\prime}$ vertices. Hence we may assume that $|(U\bigcup T)\bigcap V_2|< (1+4\xi)k^{\prime}$.
Since $|(U\bigcup S)\bigcap V_2|> (1+4\xi)k^{\prime}$, then $|S\bigcap V_2|> |T\bigcap V_2|$. Recall that $|C|> \frac{5}{2}(1+4\xi)k^{\prime}$. Thus
$$|T\bigcap V_1|> |C|-(|V_1\bigcap (U\bigcup S)|+| U\bigcap V_2|+2|S\bigcap V_2|)\geq |C|-(|U|+2|S|)> \frac{1}{2}(1+4\xi)k^{\prime}.$$
By Claim~\ref{4-10},  $|G_{B}[C\bigcap V_1,V_2\backslash S]|$ contains a connected matching saturating $(2+8\xi)k^{\prime}$ vertices. Hence we may assume that $|(T\bigcup U)\bigcap
V_1|<(1+4\xi)k^{\prime}$ and $|T\bigcap V_1|<\frac{1}{2}(1+4\xi)k^{\prime}$.

Since $|C\bigcap V_1|> \frac{3}{2}(1+4\xi)k^{\prime}$, $|S\bigcap V_1|= |C\bigcap V_1|-|(T\bigcup U)\bigcap V_1|> \frac{1}{2}(1+4\xi)k^{\prime}>|T\bigcap V_1|$. Thus
$$|T\bigcap V_2|> |C|-(|V_2\bigcap (U\bigcup S)|+| U\bigcap V_1|+2|S\bigcap V_1|)\geq |C|-(|U|+2|S|)> \frac{1}{2}(1+4\xi)k^{\prime}.$$
By Claim~\ref{4-10}, $|E(G_{B}[C\bigcap V_1,V_2\backslash S])|$ contains a connected matching saturating $(2+8\xi)k^{\prime}$ vertices. \q

\section{Remark}

Let $U$ and $V$ be vertex sets. Denote by $U\stackrel{R}{\sim} V $ (or $U\stackrel{B}{\sim} V $) be the complete bipartite graph whose edges are all colored in the red color~(or in
the blue color). Let $U=\bigcup^{4}_{i=1}U_i$, $V=\bigcup^{4}_{i=1}V_i$, and $|U_i|=|V_i|=n, i=1,2,3,4$. Take a vertex $u\in U_3$, $v\in V_1$. Now construct a $2$-colored edge
bipartite graph $\widetilde{H}$ with bipartition $\{U,V\}$ as follows:
\begin{itemize}
  \item  $U_1\stackrel{B}{\sim} V_1\bigcup V_2$, $U_1\stackrel{R}{\sim} V_4$;
  \item $U_2\stackrel{R}{\sim} V_1\bigcup V_3$, $U_2\stackrel{B}{\sim} V_2$;
  \item $U_3\stackrel{R}{\sim} V_2$, $U_3\stackrel{B}{\sim} V_3$, $U_3\backslash\{u\}\stackrel{B}{\sim} V_4$, $\{u\}\stackrel{R}{\sim} V_4$;
  \item $U_4\stackrel{R}{\sim} V_3\bigcup V_1\backslash\{v\}$, $U_4\stackrel{B}{\sim} V_4\bigcup\{v\}$.
\end{itemize}

It is clear that the bipartite graph $\widetilde{H}$ with bipartition $\{U,V\}$, $|U|=|V|=N=4n$, and $\delta(\widetilde{H})= \frac{3}{4}N$ contains neither red cycle
$C_{4n}$ nor blue cycle $C_{4n}$. We propose the following problem.

\begin{problem}\label{4-0}
Determine the smallest constant $c>0$ such that the following holds: Let $G$ be a bipartite graph with bipartition $\{V_1,V_2\}$, $|V_1|=|V_2|=N=(2+o(1))n$,
and with $\delta(G)\geq cN$ such that any $2$-edge coloring of $G$ contains a monochromatic copy of $C_{2n}$?
\end{problem}

Theorem \ref{1-5} and the above example  show that ${3 \over 4}\le c\le {7 \over 8}$. Let us ask whether  $c={3 \over 4}$ holds?


\addcontentsline{toc}{chapter}{Bibliography}
\begin{thebibliography}{99}

\bibitem{FS09} F. S. Benevides, J. Skokan, The 3-colored Ramsey number of even cycles, {\it J. Combin. Theory Ser. B,} (2009), 690-708.

\bibitem{JB73} J. A. Bondy, P. Erd\H{o}s, Ramsey numbers for cycles in graphs, {\it J. Combin. Theory Ser. B,} {\bf 14}(1973), 46-54.

\bibitem{MB18} M. Buci\'c, S. Letzter, B. Sudakov, Three colour bipartite Ramsey number of cycles and paths, arXiv:1803.03689v1.

\bibitem{CFS} D. Conlon, J. Fox and B. Sudakov, Recent developments in graph Ramsey theory. Surveys in combinatorics 2015, Cambridge University Press, pp49-118.

\bibitem{ED17} E. Davies, M. Jenssen, B. Roberts, Multicolour Ramsey numbers of paths and even cycles, {\it arXiv}:1606.00762v3, 2017.

\bibitem{LD16} L. Debiasio, L. Nelsen, Monochromatic cycle partitions of graphs with large minimum degree,  {\it J. Combin. Theory Ser. B,} {\bf 122}(2017), 634-667.

\bibitem{FT07} A. Figaj, T. {\L}uczak, The Ramsey number for a triple of long even cycles, {\it J. Combin. Theory Ser. B,} {\bf 97}(2007), 584-596.

\bibitem{RF74} R. J. Faudree, R. H. Schelp, All Ramsey numbers for cycles in graphs, {\it Discrete Math.,} {\bf 8}(1974), 313-329.

\bibitem{WG08} W. Goddard, M. A. Henning, O. R. Oellermann, Bipartite Ramsey numbers and Zarankiewicz number, {\it Discrete Math.,} {\bf 219(1)}(2008), 85-95.

\bibitem{GRS} R. Graham, B. Rothchild and J. Spencer, Ramsey theory, Wiley, New York, 1980.

\bibitem{JH18} J. H. Hattingh, E. J. Joubert, Some multicolor bipartite Ramsey numbers involving cycles and a small number of colors, {\it Discrete Math.,} {\bf 341}(2018),
    1325-1330.

\bibitem{EJ17}E. J. Joubert, Some generalized bipartite Ramsey numbers involving short cycles, {\it Graphs and Combin.,} {\bf 33} (2017), 433-448.

\bibitem{MJ16}M. Jenssen, J. Skokan, Exact Ramsey numbers of odd cycles via nonlinear optimisation, arXiv:1608.05705 (2016).

\bibitem{CK18}C. Knierim, P. Su, Improved bounds on the multicolor Ramsey numbers of paths and even cycles, arXiv:1801.04128v1(2018).

\bibitem{JK96} J. Kom\'os, M. Simonovits. Szemer\'edi's Regularity Lemma and its applications in graph thoery. {\it Combinatorics, Paul Erd\H{o}s is eighty, Vol. 2 (Keszthely,
    1993)}, 295-352, Bolyai Soc. Math. Stud., 2, {\it J\'anos Bolyai Math. Soc., Budapest,} 1996.


\bibitem{RS12} H. Li, V. Nikiforov and R. H. Schelp, A new type of  Ramsey-Tur\'an  problems, {\it Discrete Math.,} {\bf 310}(2010), 3579-3583.


\bibitem{HL09} H. Liu, Y. Person, Highly connected coloured subgraphs via the regularity lemma, {\it Discrete Math.,} {\bf 309}(2009), 6277-6287.

\bibitem{QL15}Q. Lin, Y. Li, A Folkman linear family, {\it SIAM J Discrete Math.,} {\bf 29} (2015), 1988-1998.

\bibitem{LP18} S. Liu, Y. Peng, An upper bound of bipartite Ramsey numbers of large cycles in multicolorings, manuscript.

\bibitem{TL99} T. {\L}uczak, $R(C_n,C_n,C_n)\leq(4+o(1))n$, {\it J. Combin. Theory Ser. B,} {\bf 75} (1999), 174-187.

\bibitem{VR731} V. Rosta, On a Ramsey-type problem of J. A. Bondy and P. Erd\H{o}s. I, {\it J. Combin. Theory Ser. B,} {\bf 15} (1973), 94-104.

\bibitem{VR732} V. Rosta, On a Ramsey-type problem of J. A. Bondy and P. Erd\H{o}s. II, {\it J. Combin. Theory Ser. B,} {\bf 15} (1973), 105-120.

\bibitem{GS16} G. N. S\'ark\"ozy, On the multi-colored Ramsey numbers of paths and even cycles, {\it Electron J. Combin.,} {\bf 23(3)}(2016), \#P3.53.


\bibitem{LS18} L. Shen, Q. Lin, Q. Liu, Bipartite Ramsey numbers for bipartite graphs of small bandwidth, {\it Electron J. Combin.,} {\bf 25(2)} (2018), \#P2.16.

\bibitem{RZ11} R. Zhang, Y. Sun, The bipartite Ramsey numbers $b(C_{2m};K_{2,2})$, {\it Electron J. Combin.,} {\bf 18}(2011), \#P51.

\bibitem{RZ13} R. Zhang, Y. Sun, Y. Wu, The bipartite Ramsey numbers $b(C_{2m};C_{2n})$, {\it Int. J. Math. Comp. Sci. Eng.,} {\bf 1}(2013), 80-83.





\end{thebibliography}
\end{document}